\documentclass[11pt,a4paper]{article}
\usepackage{graphicx}
\usepackage{amsmath}
\usepackage{graphicx}
\usepackage{verbatim}
\usepackage{color}
\usepackage{subfigure}
\usepackage{float}
\usepackage{authblk}
\usepackage{doi}
\usepackage{textcomp}
\usepackage{eurosym}
\newtheorem{theorem}{Theorem}[section]

\newtheorem{conjecture}{Conjecture}[section]

\begin{document}
\title{On edge irregularity strength of cycle-star graphs}
\author[1]{Umme Salma}
\author[1]{H. M. Nagesh}
\author[2]{Narahari. N}

\affil[1]{Department of Science and Humanities,\newline PES University, Bangalore, India. \newline  E-mail: ummesalma@pesu.pes.edu;  nageshhm@pes.edu \newline } 
\affil[2]{Department of Mathematics,\newline University College of Science, \newline Tumkur University, Tumakuru, India. \newline  
E-mail: narahari\_nittur@yahoo.com \newline }

\date{}

\maketitle
\begin{abstract} For a simple graph $G$, a vertex labeling $\phi:V(G) \rightarrow \{1, 2,\ldots,k\}$ is called $k$-labeling. The weight of an edge $uv$ in $G$, written $w_{\phi}(uv)$, is the sum of the labels of end vertices $u$ and $v$, i.e., $w_{\phi}(uv)=\phi(u)+\phi(v)$. A vertex $k$-labeling is defined to be an edge irregular $k$-labeling of the graph $G$ if for every two distinct edges $u$ and $v$, $w_{\phi}(u) \neq w_{\phi}(v)$. The minimum $k$ for which the graph $G$ has an edge irregular $k$-labeling is called the edge irregularity strength of $G$, written $es(G)$. In this paper, we study the edge irregular $k$-labeling for cycle-star graph $CS_{k,n-k}$ and determine the exact value for cycle-star graph for $3 \leq k \leq 7$ and $n-k \geq 1$. Finally, we make a conjecture for the edge irregularity strength of $CS_{k,n-k}$ for $k \geq 8$ and $n-k \geq 1$.  

\vskip1em \noindent \textbf{Keywords and phrases:} Irregular assignment, irregularity strength, irregular total k-labeling, edge irregularity strength, cycle-star graph.  

\vskip1em \noindent \textbf{2010 AMS Classification:} 05C38, 05C78.
\end{abstract}

\section{Introduction} \label{sec:Intr}
Let $G$ be a connected, simple, and undirected graph with vertex set $V(G)$ and edge set $E(G)$. By a labeling we mean any mapping
that maps a set of graph elements to a set of numbers (usually positive integers), called $labels$. If the domain is the vertex-
set (the edge-set), then the labeling is called \emph{vertex labelings} (\emph{edge labelings}). If the domain is $V(G) \cup E(G)$, then the labeling is called \emph{total labeling}. Thus, for an edge k-labeling $\delta: E(G) \rightarrow \{1, 2,\ldots,k\}$ the associated weight of a vertex $x \in V(G)$ is $w_{\delta}(x)=\sum \delta(xy)$, where the sum is over all vertices $y$ adjacent to $x$.
\newpage
Chartrand et al. \cite{9} defined irregular labeling for a graph $G$ as an assignment of labels from the set of natural numbers to the edges of $G$ such that the sum of the labels assigned to the edges of each vertex are different. The minimum value of the largest label of an edge over all existing irregular labelings is known as the \emph{irregularity strength} of $G$, denoted by $s(G)$. Finding the irregularity strength of a graph seems to be hard even for simple graphs, see e.g., \cite{9, 10}. 

Motivated by this, Baca et al. \cite{7} investigated the irregularity strength of graphs, namely total edge irregularity strength, denoted by $tes(G)$; and total vertex irregularity strength, denoted by $tvs(G)$. Some results on the total edge irregularity strength and the total vertex irregularity strength can be found in \cite{1, 2, 4, 6, 8}. 

Motivated by the work of Chartrand et al. \cite{9}, Ahmad et al. \cite{3} introduced the concept of \emph{edge irregular k-labelings} of graphs.

A vertex labeling $\phi:V(G) \rightarrow \{1, 2,\ldots,k\}$ is called \emph{$k$-labeling}. The weight of an edge $uv$ in $G$, written $w_{\phi}(uv)$, is the sum of the labels of end vertices $u$ and $v$, i.e., $w_{\phi}(uv)=\phi(u)+\phi(v)$. A vertex $k$-labeling of a graph $G$ is defined to be an \emph{edge irregular $k$-labeling} of the graph $G$ if for every two different edges $u$ and $v$, $w_{\phi}(u) \neq w_{\phi}(v)$. The minimum $k$ for which the graph $G$ has an edge irregular $k$-labeling is called the \emph{edge irregularity strength} of $G$, written $es(G)$. Over the last years, $es(G)$ has been investigated for different families of graphs including trees with the help of algorithmic solutions, see \cite{12, 13, 14, 16, 17, 18, 19, 20, 21, 22, 23}. The most complete recent survey of graph labelings is \cite{11}.

Sedlar \cite{15} introduced the concept of cycle-star graph as follows. The \emph{cycle-star} graph, written $CS_{k,n-k}$, is a graph with $n$ vertices consisting of the cycle graph of length $k$ and $n-k$ leafs appended to the same vertex of the cycle. 

Figure 1 shows an example of cycle-star graphs $CS_{3,4}$ and $CS_{4,3}$. 

\begin{center}
\includegraphics[width=10cm]{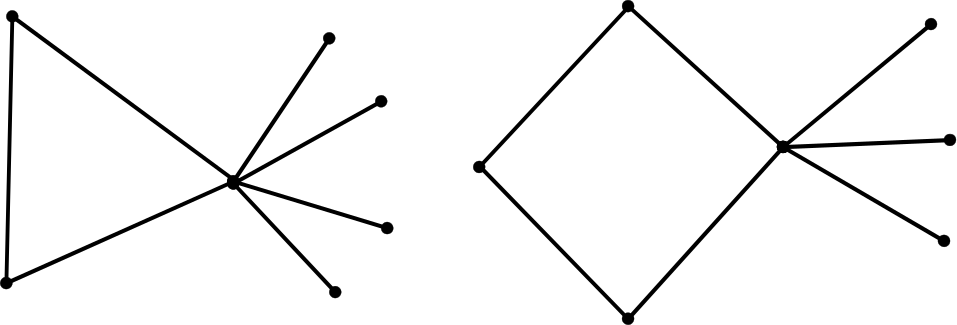}
\end{center}
\vspace{5mm}
\begin{center}
Figure 1: The cycle-star graphs $CS_{3,4}$ and $CS_{4,3}$.
\end{center}
Let the vertex set and edge set of a cycle-star graph be
\begin{center} 
$V(CS_{k,n-k})=\{v,v_i: 1 \leq i \leq n-1 \}$; 
\end{center}
\begin{equation*}
\begin{split}
E(CS_{k,n-k})  & = \{vv_i: 1 \leq i \leq n-k \} \cup \{vv_{n-k+1},vv_{n-1}\} \cup \{v_iv_{i+1}: n-k+1 \leq i \leq n-2 \}.
\end{split}
\end{equation*} 
\newpage
\section{Preliminary results}
The authors in \cite{3} estimated the bounds of the edge irregularity strength and then determined its exact values for several families
of graphs namely, paths, stars, double stars, and Cartesian product of two paths. Ahmad et al. \cite{5} determined the edge irregularity strength
of Toeplitz graphs. Recently, Tarawneh et al. \cite{18} determined the exact value of edge irregularity strength of corona product of graphs with cycle. 

Motivated by the studies mentioned above, we determine the exact value of edge irregularity strength of some classes of cycle-star graphs. 
\section{Main results}
The following theorem in \cite{3} determines the lower bound for the edge irregularity strength of a graph $G$.
\begin{theorem} Let $G=(V,E)$ be a simple graph with maximum degree $\Delta(G)$. Then \begin{center} $es(G) \geq max \{\lceil \frac{|E(G)|+1}{2} \rceil, \Delta(G)\}$
\end{center}
\end{theorem}

For two vertices $u$ and $v$ in a graph $G$, the $distance$ $d(u,v)$ from $u$ to $v$ is the length of a shortest $u-v$ path in $G$. For a vertex $v$ in a connected graph $G$, the $eccentricity$ $e(v)$ of $v$ is the distance between $v$ and a vertex farthest from $v$ in $G$. The minimum eccentricity among the vertices of $G$ is its $radius$ and the maximum eccentricity is its $diameter$, which are denoted by $rad(G)$ and $diam(G)$, respectively. A vertex $v$ in $G$ is a \emph{central vertex} if $e(v)=rad(G)$.

We will make use of Theorem 3.1 to prove our main results. In the next theorem, we determine the exact value of the edge irregularity strength of cycle-star graph $CS_{k,n-k}$ for $k=3$ and $n-k \geq 1$.
\begin{theorem}
Let $G=CS_{k,n-k}$ be a cycle-star graph. For $k=3$ and $n-k \geq 1$, $es(G)=n-1$.
\end{theorem}
\textbf{Proof}. Let $G=CS_{k,n-k}$ be a cycle-star graph, where $k=3$ and $n-k \geq 1$. Let us consider the vertex set and edge set of $G$. 
\begin{center} 
$V(G)=\{v,v_i: 1 \leq i \leq n-1 \}$; 
\end{center}
\begin{equation*}
\begin{split}
E(G)  & = \{vv_i: 1 \leq i \leq n-3 \} \cup \{vv_{n-2},vv_{n-1}\} 
 \cup \{v_{n-2}v_{n-1}\}.
\end{split}
\end{equation*}
Let $v$ be the central vertex; $v_1,v_2,\ldots,v_{n-3}$ be leafs that are adjacent to $v$; and let $v_{n-1},v_{n-2}$ be the other vertices on the cycle $C_3$. We have to show that $es(G)=n-1$. From Theorem 3.1, we get the lower bound $es(G)\geq n-1$.

To prove the equality, it suffices to prove the existence of an edge irregular $(n-1)-$ labeling. 

Define a labeling on vertex set of $G$ as follows:\\
Let $\phi:V(G) \rightarrow \{1,2,\ldots,n-1\}$ such that 
$\phi(v)=1$ and $\phi(v_i)=i$ for $i=1,2,\ldots,n-1$. \newpage The edge weights are as follows:
\begin{center}
$w_{\phi}(vv_i)=1+i  \hspace{3mm} \text{for }  i=1,2,\ldots,n-3$.
\end{center}
\begin{center}
$w_{\phi}(vv_{n-2})=n-1; w_{\phi}(vv_{n-1})=n; w_{\phi}(v_{n-2}v_{n-1})=2n-3$.
\end{center}
On the basis of above calculations we see that the edge weights are distinct for all pairs of distinct edges. Thus the vertex labeling $\phi$ is an edge regular $(n-1)-$ labeling. Therefore, $es(G)=n-1$.

\begin{center}
\includegraphics[width=85mm]{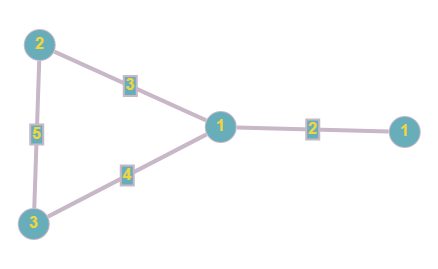}
\end{center}
\begin{center}
Figure 2: Edge irregularity strength of $CS_{3,1}$.
\end{center}

In the next theorem, we determine the exact value of the edge irregularity strength of $CS_{k,n-k}$ for $k=4$ and $n-k \geq 1$.
\begin{theorem}
Let $G=CS_{k,n-k}$ be a cycle-star graph. For $k=4$ and $n-k \geq 1$, $es(G)=n-2$.
\end{theorem}
\textbf{Proof}.
Let $G=CS_{k,n-k}$ be a cycle-star graph, where $k=4$ and $n-k \geq 1$. 
Let us consider the vertex set and edge set of $G$. 
\begin{center} 
$V(G)=\{v,v_i: 1 \leq i \leq n-1 \}$; 
\end{center}
\begin{equation*}
\begin{split}
E(G)  & = \{vv_i: 1 \leq i \leq n-4 \} \cup \{vv_{n-3},vv_{n-1}\} 
 \cup \{v_{i}v_{i+1}: n-3 \leq i \leq n-2\}.
\end{split}
\end{equation*}
Let $v$ be the central vertex; $v_1,v_2,\ldots,v_{n-4}$ be leafs that are adjacent to $v$; and let $v_{n-1},v_{n-2},v_{n-3}$ be the other vertices on the cycle $C_4$. We have to show that $es(G)=n-2$. From Theorem 3.1, we get the lower bound $es(G)\geq n-2$.

To prove the equality, it suffices to prove the existence of an edge irregular $(n-2)-$ labeling. 

Define a labeling on vertex set of $G$ as follows:\\
Let $\phi:V(G) \rightarrow \{1,2,\ldots,n-2\}$ be such that $\phi(v)=1$; $\phi(v_i)=i$ for $i=1,2,\ldots,n-2$; and $\phi(v_{n-1})=3$. \newpage The edge weights are as follows:
\begin{center}
$w_{\phi}(vv_i)=i+1  \hspace{3mm} \text{for }  i=1,2,\ldots,n-4
$.
\end{center}
\begin{center}
$w_{\phi}(vv_{n-3})=n-2; w_{\phi}(vv_{n-2})=n-1$.
\end{center}
\begin{center}
$w_{\phi}(v_{n-3}v_{n-1})=n; w_{\phi}(v_{n-2}v_{n-1})=n+1$.
\end{center}
On the basis of above calculations we see that the edge weights are distinct for all pairs of distinct edges. Thus the vertex labeling $\phi$ is an edge regular $(n-2)-$ labeling. Therefore, $es(G)=n-2$.

\begin{center}
\includegraphics[width=100mm]{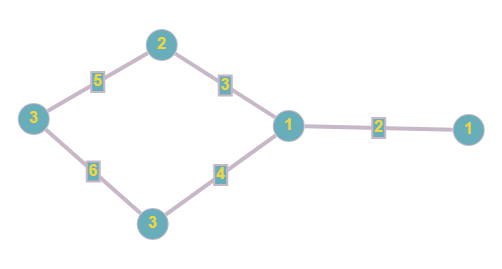}
\end{center}
\begin{center}
Figure 3: Edge irregularity strength of $CS_{4,1}$.
\end{center}
In the next theorem, we determine the exact value of the edge 
irregularity strength of cycle-star graph $CS_{k,n-k}$ for $k=5$ and $n-k \geq 1$.
\begin{theorem}
Let $G=CS_{k,n-k}$ be a cycle-star graph. For $k=5$ and $n-k \geq 1$,
\begin{center}
$es(G):=
\begin{cases}
n-2 & \text{for $n=6$  } 
\\
n-3 & \text{for $n \geq 7$} 
\end{cases}$
\end{center} 
\end{theorem}
\textbf{Proof}.
Let $G=CS_{k,n-k}$ be a cycle-star graph, where $k=5$ and $n-k \geq 1$. 
Let us consider the vertex set and edge set of $G$. 
\begin{center} 
$V(G)=\{v,v_i: 1 \leq i \leq n-1 \}$; 
\end{center}
\begin{equation*}
\begin{split}
E(G)  & = \{vv_i: 1 \leq i \leq n-5 \} \cup \{vv_{n-4},vv_{n-1}\} 
 \cup \{v_{i}v_{i+1}: n-4 \leq i \leq n-2\}.
\end{split}
\end{equation*}
We consider the following two cases.\\
\textbf{Case $1$}: The graph $G=CS_{5,n-5}$, where $n \geq 6$, is of order $n$, size $n$, and has the maximum $\Delta=n-3$. Thus by Theorem 3.1, we get $es(G) \geq max \{\lceil \frac{n+1}{2} \rceil, n-3\}$. 

For $n=6$, $n-3=6-3=3<\lceil \frac{n+1}{2} \rceil=\lceil \frac{6+1}{2} \rceil=4$. Thus $es(G)  \geq 4 \, (=n-2)$.

To prove the equality, it suffices to prove the existence of an edge irregular $(n-2)-$ labeling.

For the cycle-star graph $CS_{5,1}$, let $v$ be the central vertex; $v_1$ be the leaf adjacent to $v$; and let $v_{2},v_{3},v_{4},v_{5}$ be the other vertices on cycle $C_5$.

Let $\phi:V(CS_{5,1}) \rightarrow \{1,2,3,4\}$ be the vertex labeling such that $\phi(v)=1$; $\phi(v_i)=i$ for $1 \leq i \leq 4 $; and $\phi(v_{5})=4$.

The edge weights are as follows:
\begin{center}
$w_{\phi}(vv_i)=i+1  \hspace{3mm} \text{for }  i=1,2,3$.
\end{center}
\begin{center}
$w_{\phi}(v_{2}v_{4})=6; w_{\phi}(v_{3}v_{5})=7; w_{\phi}(v_{4}v_{5})=8$.
\end{center}
On the basis of above calculations we see that the edge weights are distinct for all pairs of distinct edges. Thus the vertex labeling $\phi$ is an edge regular $(n-2)-$ labeling. Therefore, $es(CS_{5,1})=n-2$ for $n=6$.
\vspace{5mm}
\begin{center}
\includegraphics[width=75mm]{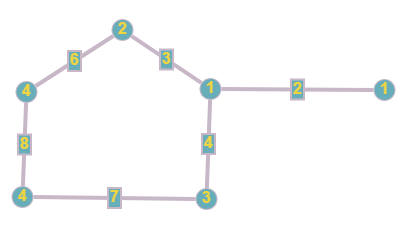}
\end{center}
\begin{center}
Figure 4: Edge irregularity strength of $CS_{5,1}$.
\end{center}
\textbf{Case $2$}: For the cycle-star graph $G=CS_{5,n-5}$, where $n \geq 7$, it is easy to see that $n-3 \geq \lceil \frac{n+1}{2} \rceil$. Thus, by Theorem 3.1, $es(G)  \geq n-3$. To prove the equality, it suffices to prove the existence of an edge irregular $(n-3)-$ labeling.

For the cycle-star graph $G=CS_{5,n-5}$, where $n \geq 7$, let $v$ be the central vertex; $v_1,v_2,\ldots,v_{n-5}$ be leafs that are adjacent to $v$; and let $v_{n-1},v_{n-2},v_{n-3},v_{n-4}$ be the other vertices on the cycle $C_5$. 

Let $\phi:V(G) \rightarrow \{1,2,\ldots,n-3\}$ be the vertex labeling such that $\phi(v)=1$; $\phi(v_i)=i$ for $i=1,2,\ldots,n-3$; $\phi(v_{n-2})=n-4$; and $\phi(v_{n-1})=n-3$. 

The edge weights are as follows:
\begin{center}
$w_{\phi}(vv_i)=i+1  \hspace{3mm} \text{for }  i=1,2,\ldots,n-3$.
\end{center}
\begin{center}
$w_{\phi}(v_{n-4}v_{n-2})=2(n-4); w_{\phi}(v_{n-2}v_{n-1})=2n-7$.
\end{center}
\begin{center}
$w_{\phi}(v_{n-3}v_{n-1})=2(n-3)$.
\end{center}
On the basis of above calculations we see that the edge weights are distinct for all pairs of distinct edges. Thus the vertex labeling 
$\phi$ is an edge regular $(n-3)-$ labeling. Therefore, 
$es(G)=n-3$ for $k=5$ and $n \geq 7$.

\vspace{5mm}
\begin{center}
\includegraphics[width=100mm]{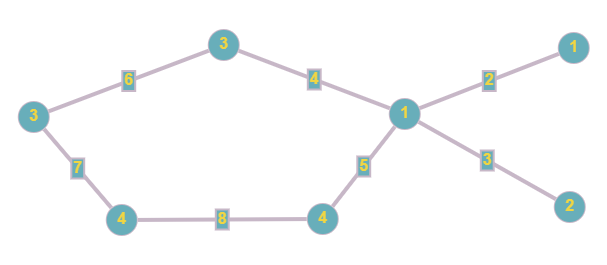}
\end{center}
\begin{center}
Figure 5: Edge irregularity strength of $CS_{5,2}$.
\end{center}
In the next theorem, we determine the exact value of the edge irregularity strength of cycle-star graph $CS_{k,n-k}$ for $k=6$ and $n-k \geq 1$. 
\begin{theorem}
Let $G=CS_{k,n-k}$ be a cycle-star graph. For $k=6$ and $n-k \geq 1$, 
\begin{center}
$es(G):=
\begin{cases}
n-3 & \text{for $n=7,8$  } 
\\
n-4 & \text{for $n \geq 9$} 
\end{cases}$
\end{center} 
\end{theorem}
\textbf{Proof}.
Let $G=CS_{k,n-k}$ be a cycle-star graph, where $k=6$ and $n-k \geq 1$. Let us consider the vertex set and edge set of $G$. 
\begin{center} 
$V(G)=\{v,v_i: 1 \leq i \leq n-1 \}$; 
\end{center}
\begin{equation*}
\begin{split}
E(G)  & = \{vv_i: 1 \leq i \leq n-6 \} \cup \{vv_{n-5},vv_{n-1}\} 
 \cup \{v_{i}v_{i+1}: n-5 \leq i \leq n-2\}.
\end{split}
\end{equation*} 
We consider the following two cases.\\
\textbf{Case $1$}: The cycle-star graph $G=CS_{6,n-6}$, where $n \geq 7$, is of order $n$, size $n$, and has the maximum $\Delta=n-4$. Thus by Theorem 3.1, we get $es(G) \geq max \{\lceil \frac{n+1}{2} \rceil, n-4\}$. 

For $n=7$, $n-4=7-4=3<\lceil \frac{n+1}{2} \rceil=\lceil \frac{7+1}{2} \rceil=4$. Thus $es(G)  \geq 4 \, (=n-3)$. Again, for $n=8$, $n-4=8-4=4<\lceil \frac{n+1}{2} \rceil=\lceil \frac{8+1}{2} \rceil=5$. Thus $es(G)  \geq 5 \, (=n-3)$. Therefore, in both the cases, $es(CS_{6,n-6}) \geq n-3$ for $n=7,8$.

To prove the equality, it suffices to prove the existence of an edge irregular $(n-3)-$ labeling.

For the cycle-star graph $G=CS_{6,n-6}$, where $n=7,8$, let $v$ be the central vertex; $v_2$ be the leaf adjacent to $v$; and let $v_{1},v_{3},v_{4},v_{5},\ldots,v_{n-1}$ be the other vertices on the cycle $C_6$.

Let $\phi:V(G) \rightarrow \{1,2,\ldots,n-3\}$ be the vertex labeling such that $\phi(v)=1$; $\phi(v_i)=i$ for $1 \leq i \leq n-4 $; $\phi(v_{n-3})=\phi(v_{n-2})=n-3$; and $\phi(v_{n-1})=n-4$. 

The edge weights are as follows:
\begin{center}
$w_{\phi}(vv_i)=i+1  \hspace{3mm} \text{for }  2 \leq i \leq n-4$;
\end{center}
\begin{center}
$w_{\phi}(v_{1}v_{n-3})=n-2; w_{\phi}(v_{n-3}v_{n-2})=2(n-3)$;
\end{center}
\begin{center}
$w_{\phi}(v_{n-2}v_{n-1})=2n-7; w_{\phi}(v_{n-4}v_{n-1})=2(n-4)$.
\end{center}
On the basis of above calculations we see that the edge weights are distinct for all pairs of distinct edges. Thus the vertex labeling 
$\phi$ is an edge regular $(n-3)-$ labeling. Therefore, $es(G)=n-3$ for $n=7,8$.
\begin{center}
\includegraphics[width=150mm]{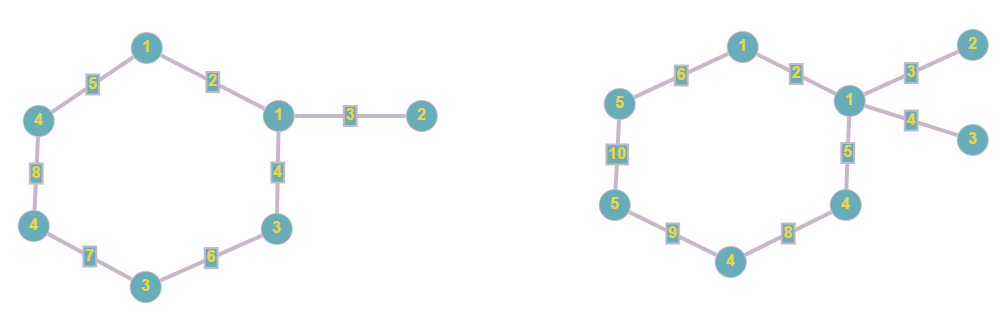}
\end{center}
\begin{center}
Figure 6: Edge irregularity strength of $CS_{6,1}$ and $CS_{6,2}$.
\end{center}
\textbf{Case $2$}: For the cycle-star graph $G=CS_{6,n-6}$, where $n \geq 9$, it is easy to see that $n-4 \geq \lceil \frac{n+1}{2} \rceil$. Thus, by Theorem 3.1, $es(G)  \geq n-4$.

To prove the equality, it suffices to prove the existence of an edge irregular $(n-4)-$ labeling.

For the cycle-star graph $G=CS_{6,n-6}$, where $n \geq 9$, let $v$ be the central vertex; $v_2,v_3,\ldots,v_{n-5}$ be leafs that are adjacent to $v$; and let $v_{1},v_{n-4},v_{n-3},v_{n-2},v_{n-1}$ be the other vertices on the cycle $C_6$.

Let $\phi:V(G) \rightarrow \{1,2,\ldots,n-4\}$ be the vertex labeling such that $\phi(v)=n-4$; $\phi(v_i)=i$ for $1 \leq i \leq 2 $; $\phi(v_i)=i+1$ for $3 \leq i \leq n-5 $; $\phi(v_{n-4})=3$; $\phi(v_{n-3})=\phi(v_{n-2})=2$; and $\phi(v_{n-1})=1$. 

The edge weights are as follows:
\begin{center}
$w_{\phi}(vv_i)=i+n-4  \hspace{3mm} \text{for } 1 \leq i \leq 2$. \\
\end{center}
\begin{center}
$w_{\phi}(vv_{i})=n-1+j \hspace{3mm} \text{for } 3 \leq i \leq n-5 \hspace{2mm} \text{and } 1 \leq j \leq n-7$.
\end{center}
\begin{center}
$w_{\phi}(vv_{n-4})=n-1; w_{\phi}(v_{1}v_{n-1})=2$.
\end{center}
\begin{center}
$w_{\phi}(v_{n-4}v_{n-3})=5; w_{\phi}(v_{n-3}v_{n-2})=4; w_{\phi}(v_{n-2}v_{n-1})=3$.
\end{center}
On the basis of above calculations we see that the edge weights are distinct for all pairs of distinct edges. Thus the vertex labeling 
$\phi$ is an edge regular $(n-4)-$ labeling. Therefore, $es(G)=n-4$ for $k=6$ and $n \geq 9$.

\begin{center}
\includegraphics[width=100mm]{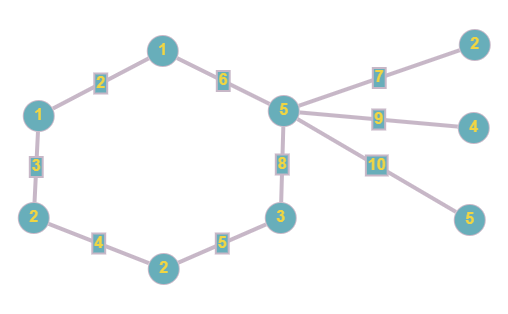}
\end{center}
\begin{center}
Figure 7: Edge irregularity strength of $CS_{6,3}$.    
\end{center}
In the next theorem, we determine the exact value of the edge irregularity strength of cycle-star graph $CS_{k,n-k}$ for $k=7$ and $n-k \geq 1$.
\begin{theorem}
Let $G=CS_{k,n-k}$ be a cycle-star graph. For $k=7$ and $n-k \geq 1$, 
\begin{center}
$es(G):=
\begin{cases}
n-3 & \text{for $n=8$  } 
\\
n-4 & \text{for $n=9,10$  } 
\\
n-5 & \text{for $n \geq 11$} 
\end{cases}$
\end{center} 
\end{theorem}
\textbf{Proof}.
Let $G=CS_{k,n-k}$ be a cycle-star graph, where $k=7$ and $n-k \geq 1$. 
Let us consider the vertex set and edge set of $G$. 
\begin{center} 
$V(G)=\{v,v_i: 1 \leq i \leq n-1 \}$; 
\end{center}
\begin{equation*}
\begin{split}
E(G)  & = \{vv_i: 1 \leq i \leq n-6 \} \cup \{vv_{n-6},vv_{n-1}\} 
 \cup \{v_{i}v_{i+1}: n-6 \leq i \leq n-2\}.
\end{split}
\end{equation*} 
We consider the following three cases.\\
\textbf{Case $1$}: The cycle-star graph $G=CS_{7,n-7}$, where $n \geq 8$, is of order $n$, size $n$, and has the maximum $\Delta=n-5$. Thus by Theorem 3.1, we get $es(G) \geq max \{\lceil \frac{n+1}{2} \rceil, n-5\}$. 

For $n=8$, $n-5=8-5=3<\lceil \frac{n+1}{2} \rceil=\lceil \frac{8+1}{2} \rceil=5$. Thus $es(G)  \geq 5 \, (=n-3)$. 
\newpage
To prove the equality, it suffices to prove the existence of an edge irregular $(n-3)-$ labeling. 

For the cycle-star graph $g=CS_{7,1}$, let $v$ be the central vertex; $v_1$ be the leaf adjacent to $v$; and let  $v_{2},v_{3},v_{4},v_{5},v_{6},v_{7}$ be the other vertices on the cycle $C_7$. 

Let $\phi:V(G) \rightarrow \{1,2,\ldots,n-3\}$ be the vertex labeling such that $\phi(v)=1$; $\phi(v_i)=i$ for $ 1 \leq i \leq 3$; 
$\phi(v_{4})=\phi(v_{5})=3$; $\phi(v_{6})=4$; and $\phi(v_{7})=5$. 

The edge weights are as follows:
\begin{center}
$w_{\phi}(vv_i)=i+1  \hspace{3mm} \text{for }  1 \leq i \leq 3$.
\end{center}
\begin{center}
$w_{\phi}(v_{2}v_{4})=5; w_{\phi}(v_{4}v_{5})=6; w_{\phi}(v_{5}v_{6})=7; w_{\phi}(v_{6}v_{7})=9$; and $w_{\phi}(v_{3}v_{7})=8$.
\end{center}
On the basis of above calculations we see that the edge weights are distinct for all pairs of distinct edges. Thus the vertex labeling 
$\phi$ is an edge regular $(n-3)-$ labeling. Therefore, $es(G)=n-3$. 

\begin{center}
\includegraphics[width=100mm]{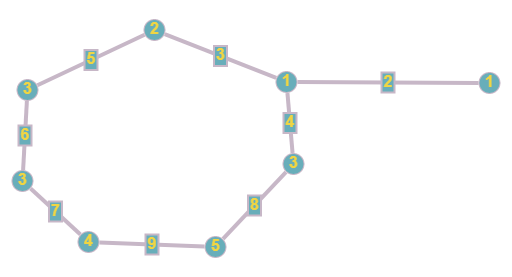}
\end{center}
\begin{center}
Figure 8: Edge irregularity strength of $CS_{7,1}$.
\end{center}
\textbf{Case $2$}: For $n=9$, $n-5=9-5=4<\lceil \frac{n+1}{2} \rceil=\lceil \frac{9+1}{2} \rceil=5$. Thus $es(CS_{7,2})  \geq 5 \, (=n-4)$. 
Again, for $n=10$, $n-5=10-5=5<\lceil \frac{n+1}{2} \rceil=\lceil \frac{10+1}{2} \rceil=6$. Thus $es(CS_{7,3})  \geq 6 \, (=n-4)$. Therefore, in both the cases, $es(G) \geq n-4$ for $n=9,10$.

To prove the equality, it suffices to prove the existence of an edge irregular $(n-4)-$ labeling.

For the cycle-star graph $G=CS_{7,n-7}$, where $n=9,10$, let $v$ be the central vertex; $v_1,v_2,\ldots,v_{n-7}$ be leafs that are adjacent to $v$; and let $v_{n-6},v_{n-5},v_{n-4},v_{n-3},v_{n-2},v_{n-1}$ be the other vertices on the cycle $C_7$.

Let $\phi:V(G) \rightarrow \{1,2,\ldots,n-4\}$ be the vertex labeling such that $\phi(v)=5$; $\phi(v_1)=3$; $\phi(v_i)=i+3$ for 
$2 \leq i \leq n-7 $; $\phi(v_{n-5})=4$; $\phi(v_{n-6})=\phi(v_{n-4})=1$; $\phi(v_{n-3})=\phi(v_{n-2})=2$; and $\phi(v_{n-1})=3$. 

The edge weights are as follows:
\begin{center}
$w_{\phi}(vv_1)=8; w_{\phi}(vv_i)=8+i  \hspace{3mm} \text{for } 2 \leq i \leq n-7$.
\end{center}
\begin{center}
$w_{\phi}(vv_{n-6})=6; w_{\phi}(vv_{n-5})=9; w_{\phi}(v_{n-5}v_{n-1})=7; w_{\phi}(v_{n-6}v_{n-4})=2$.
\end{center}
\begin{center}
$w_{\phi}(v_{i}v_{i+1})=2+j \hspace{5mm} \text{for } n-4 \leq i \leq n-2 \hspace{2mm} \text{and } 1 \leq j \leq 3$.
\end{center}
On the basis of above calculations we see that the edge weights are distinct for all pairs of distinct edges. Thus the vertex labeling 
$\phi$ is an edge regular $(n-4)-$ labeling. Therefore, $es(G)=n-4$ for $n=9,10$.
\begin{center}
\includegraphics[width=150mm]{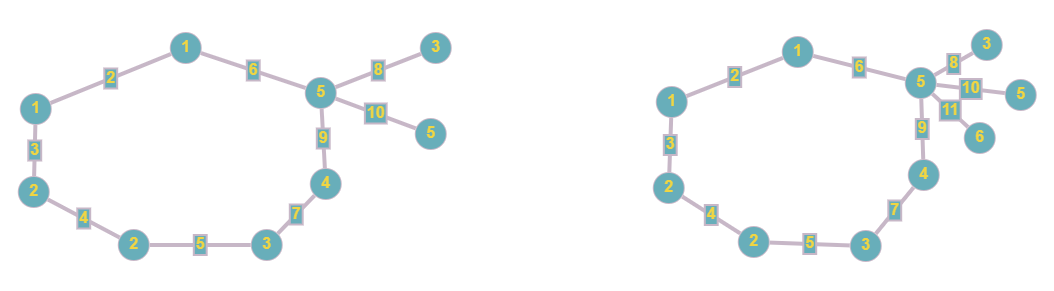}
\end{center}
\begin{center}
Figure 9: Edge irregularity strength of $CS_{7,2}$ and $CS_{7,3}$.
\end{center}
\textbf{Case $3$}: For $n \geq 11$, $n-5 \geq \lceil \frac{n+1}{2} \rceil$. Thus, by Theorem 3.1, $es(CS_{7,n-7})  \geq n-5$. 

To prove the equality, it suffices to prove the existence of an edge irregular $(n-5)-$ labeling.

For the cycle-star graph $G=CS_{7,n-7}$, where $n \geq 11$, let $v$ be the central vertex; $v_2,v_4,v_5,v_6,\ldots,v_{n-5}$ be leafs that are adjacent to $v$; and let $v_{1},v_{n-4},v_{n-3},v_{n-2},v_{n-1},v_3$ be the other vertices on the cycle $C_7$.

Let $\phi:V(G) \rightarrow \{1,2,\ldots,n-5\}$ be the vertex labeling such that $\phi(v)=n-5$; $\phi(v_i)=i$ for $1 \leq i \leq n-5$; $\phi(v_{n-4})=1$; $\phi(v_{n-3})=\phi(v_{n-2})=2$; and 
$\phi(v_{n-1})=3$. 

The edge weights are as follows:
\begin{center}
$w_{\phi}(vv_i)=n-1+i  \hspace{3mm} \text{for } 1 \leq i \leq n-5$.
\end{center}
\begin{center}
$w_{\phi}(v_{1}v_{n-4})=2; w_{\phi}(v_{n-4}v_{n-3})=3; w_{\phi}(v_{n-3}v_{n-2})=4$.
\end{center}
\begin{center}
$w_{\phi}(v_{n-2}v_{n-1})=5; w_{\phi}(v_{n-3}v_{3})=6$.
\end{center}
On the basis of above calculations we see that the edge weights are distinct for all pairs of distinct edges. Thus the vertex labeling 
$\phi$ is an edge regular $(n-5)-$ labeling. Therefore, 
$es(G)=n-5$ for $k=7$ and $n \geq 11$.

\begin{center}
\includegraphics[width=120mm]{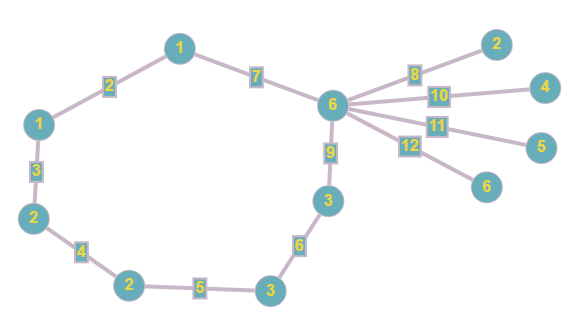}
\end{center}
\begin{center}
Figure 10: Edge irregularity strength of $CS_{7,4}$. 
\end{center}
\vspace{5mm}
We close with the following conjecture. 
\begin{conjecture}
Let $G=CS_{k,n-k}$ be a cycle-star graph. For $k \geq 8$ and $n-k \geq 1$, 
\begin{center}
$es(G):=
\begin{cases}
\lceil \frac{n+1}{2} \rceil & \text{for \, $k+1 \leq n \leq 2k-4$  } 
\\
n-k+2 & \text{for \, $n \geq 2k-3$  }  
\end{cases}$
\end{center} 
\end{conjecture}
\section{Conclusion} 
In this paper, we investigated the edge irregularity strength, as a modification of the well-known irregularity strength, total edge 
irregularity strength and total vertex irregularity strength. We obtained the exact values for edge irregularity strength of cycle-star graph $CS_{k,n-k}$ for $ 3 \leq k \leq 7$ and $n-k \geq 1$. Also, we conjectured the edge irregularity strength of cycle-star graph $CS_{k,n-k}$ for $k \geq 8$ and $n-k \geq 1$. The exact values of edge irregularity strength can be determined for graph operations, graph products, and graph powers also.

\makeatletter
\renewcommand{\@biblabel}[1]{[#1]\hfill}
\makeatother

\end{document}